\documentclass[11pt, leqno]{scrartcl}
\usepackage{amssymb,amsmath,amsthm}
\usepackage{latexsym}
\usepackage{t1enc}
\usepackage{lmodern}
\usepackage[utf8]{inputenc}
\usepackage{eurosym}
\usepackage{mdwlist}
\usepackage{paralist}
\usepackage{enumitem}
\usepackage{esint}
\usepackage{nicefrac}

\newtheorem{theorem}{Theorem}

\newtheorem{corollary}{Corollary}
\newtheorem{remark}{Remark}
\newtheorem{lemma}{Lemma}

\usepackage[pdfproducer={LaTeX with hyperref package},pdfpagelayout=SinglePage,bookmarksnumbered,pdftex=true,breaklinks=true,bookmarks=true,linktocpage=true,pdfpagelabels=true,plainpages=false,pagebackref=false]{hyperref} 
\hypersetup{colorlinks=true,linkcolor=blue,citecolor=blue,urlcolor=blue,filecolor=blue}

\newcommand{\eps}{\varepsilon}
\newcommand{\Ge}{\Gamma_{\eps,k}}
\newcommand{\Gei}[1]{\Gamma_{\eps,#1}}
\newcommand{\G}{\Gamma} 
\newcommand{\Iei}[1]{\mathcal{I}_{\eps,#1}} 
\newcommand{\uei}[1]{u_{\eps,#1}}	
\newcommand{\ue}{u_{\eps,k}}	   
\newcommand{\pen}[1]{p_{\eps,#1}}	
\newcommand{\lap}[1]{\Delta \hspace{-0.15mm} #1 }
\newcommand{\dist}{\operatorname{dist}}
\newcommand{\bilap}[1]{\Delta^{\!2}\hspace{-0.25mm}#1 }
\newcommand{\grad}[1]{\nabla \! #1} 
\newcommand{\aq}{\Leftrightarrow} 
\renewcommand{\phi}{\varphi}   
\newcommand{\norm}[1]{\|#1\|} 
\newcommand{\abs}[1]{|#1|}
\newcommand{\R}{\mathbb{R}}   
\newcommand{\N}{\mathbb{N}}   
\newcommand{\Ray}{\mathcal{R}}
\renewcommand{\O}{\mathcal{O}}

\begin{document}
\title{Existence of an optimal domain for minimizing the fundamental tone of a clamped plate of prescribed volume in arbitrary dimension}
\author{Kathrin Stollenwerk} \date{}
\pagestyle{myheadings}
\maketitle

\noindent\textbf{Abstract.} In the 19th century, Lord Rayleigh conjectured that among all clamped plates with given area, the disk minimizes the fundamental tone \cite{rayleigh}. In the 1990s, N. S. Nadirashvili \cite{nadirashvili} proved the conjecture in $\R^2$ and M. S. Ashbaugh und R. D. Benguria \cite{AshBeng1995_clamped} gave a proof in $\R^2$ and $\R^3$. 

In the present paper, we prove existence of an optimal domain for minimizing the fundamental tone among all open and bounded subsets of $\R^n$, $n\geq 4$, with given measure. 
We formulate the minimization of the fundamental tone of a clamped plate as a free boundary value problem with a penalization term for the volume constraint. As the penalization parameter becomes small we show that the optimal shape problem  is solved. 

\section{Introduction}

We consider the following variational problem. For $n\geq 2$ let $\Omega \subset \R^n$ be an open and bounded domain. Then we define
\[
\Ray(v,\Omega) := \frac{\int_\Omega \abs{\lap v}^2dx}{\int_
\Omega v^2dx}
\]
for $v \in H^{2,2}_0(\Omega)$ and denote $\Ray(v,\Omega)=\infty$ if the denominator vanishes. The quantity 
\[
\G(\Omega) := \min \{ \Ray(v,\Omega) : v \in H^{2,2}_0(\Omega)\}
\]
is called the fundamental tone of $\Omega$. The infimum is attained by the first eigenfunction which solves the following Euler-Lagrange equation
\begin{equation}\label{eq:cl_ef}
 \begin{cases}
   \bilap u  - \G(\Omega)\,u &=0 \mbox{ in } \Omega \\
   u = \abs{\grad u} &= 0 \mbox{ on } \partial\Omega
   \end{cases}
 \end{equation}
 if the boundary of $\Omega$ is smooth enough. We will denote a solution $u$ of \eqref{eq:cl_ef} a clamped eigenfunction on $\Omega$. 
 
Already in the 19th century, Lord Rayleigh conjectured that among all domains of given measure the ball has the minimal fundamental tone \cite{rayleigh}. 
However, up to now, Lord Rayleigh's conjecture is still not completely proven in any dimension. 

In 1981, G. Talenti \cite{talenti} applied rearrangement techniques to prove that, in the two-dimensional case, Rayleigh's conjecture holds true up to a factor of $0.98$. Later, N. S. Nadirashvili improved Talenti's approach and proved Rayleigh's conjecture for $n=2$ \cite{nadirashvili}. Also basing on Talenti's work, M. S. Ashbaugh and R. D. Benguria gave a proof for Rayleigh's conjecture which holds true in dimension $n=2$ and $n=3$ \cite{AshBeng1995_clamped}. 

For dimensions  $n\geq 4$, there is only a partial result known. Assuming that the eigenfunction $u$ is of fixed sign, G. Szegö showed that Rayleigh's conjecture holds true \cite{Szego50}. However, this assumption does not hold true in general. For counter examples we refer to \cite[Section 8]{AshLau1996}.  

In the present paper, we will contribute to the partial answers to Rayleigh's conjecture in dimension $n\geq 4$ and will prove the following main theorem. 

\begin{theorem}\label{theo:main}
  There exists a bounded domain $\Omega^\ast$ with given measure $\abs{\Omega}=\omega_0$ such that 
  \[
   \G(\Omega^\ast) = \min\{\G(D):  D \subset B, D \mbox{ open}, \abs{D}\leq \omega_0\}.
  \]
  Thereby, $n\geq 4$, $\omega_0>0$ is a given quantity, $\abs{ D}$ denotes the $n$-dimensional Lebesgue-measure of $D$ and $B$ denotes a ball in $\R^n$ with $\abs{B} \gg \omega_0$.
\end{theorem}

In order to prove Theorem \ref{theo:main}, we will follow an idea of H. W. Alt and L. A. Caffarelli in \cite{AltCaf81} and introduce a penalized variational problem. Let $n\geq 4$ and $0<\omega_0$. Now choose a ball $B\subset \R^n$ with $\omega_0 \ll \abs{B}$, where $\abs{B}$ denotes the $n$-dimensional Lebesgue measure of $B$. For $\eps>0$ we define the penalization term $\pen{0}: \R \to\R$ by
\[
\pen{0}(s) := \begin{cases}
  0, &s\leq \omega_0 \\
  \frac{1}{\eps}(s-\omega_0), &s\geq \omega_0
\end{cases}
\]
and the functional $\Iei{0} : H^{2,2}_0(B) \to \R$ by
\[
  \Iei{0} := \Ray(v,\Omega) + \pen{0}(\abs{\O(v)}),
\]
where $\O(v) := \{x \in B: v(x)\neq 0\}$. We will see that for each $\eps>0$ there exists a function $\uei{0} \in H^{2,2}_0(B)$, which minimizes the functional $\Iei{0}$. Provided that $\eps$ is chosen sufficiently small, the minimizer $\uei{0}$ yields a domain $\Omega(\uei{0})$ with $\abs{\Omega(\uei{0})}=\omega_0$ which minimizes the fundamental tone $\G$ among all open subsets of $B$ with measure smaller or equal than $\omega_0$. This proves the main theorem, Theorem \ref{theo:main}. 

Facing Rayleigh's conjecture, the next reasonable step should be the analysis of regularity properties of the free boundary $\partial\Omega(\uei{0})$. If we orientate ourselves on the seminal work of Alt and Caffarelli in \cite{AltCaf81}, the next aim would be establishing a nondegeneracy result for the minimizing function $\uei{0}$ along $\partial\Omega(\uei{0})$ and, subsequently, proving that $\partial\Omega(\uei{0})$ has got a positive Lebesgue density in every point. However, looking at nondegeneracy results for second order problems (see \cite{AguAltCaf1986,AltCaf81,BaWa09}, e.g.),  we see that in these settings the desired nondegeneracy result is achieved by constructing suitable testfunctions which heavily rely on comparison principles. Since such comparison principles are not available for fourth order operators in general, this might be the end point for our approach via the functional $\Iei{0}$, 
   
Consequently, we revise the functional $\Iei{0}$ and replace the penalization term. Following an idea of N. Aguilera, H. W. Alt and L. A. Caffarelli in \cite{AguAltCaf1986}, we define the penalization term $\pen{1}: \R\to\R$ by 
\[
 \pen{1}(s) := \begin{cases}
   \eps(s-\omega_0), & s\leq \omega_0 \\
   \frac{1}{\eps}(s-\omega_0), & s\geq \omega_0 
 \end{cases}
\]
and $\Iei{1} : H^{2,2}_0(B) \to \R$ by 
\[
  \Iei{1}(v) := \Ray(v,B) + \pen{1}(\abs{\O(v)}).
\]
The new penalization term rewards volumes less than $\omega_0$ with a negative contribution to the functional. We will prove that for every $\eps>0$, there exists a minimizer $\uei{1} \in H^{2,2}_0(B)$ of $\Iei{1}$ and that $\uei{1}$ yields a domain $\Omega(\uei{1})$ which minimizes the fundamental tone among all open subsets of $B$ with the same measure as $\Omega(\uei{1})$. The analysis of the volume of $\Omega(\uei{1})$ will be more challenging than in the previous case since the rewarading property of the $\pen{1}$ counteracts the monotonicity of the fundamental tone with respect to set inclusion. 

The present paper is organized as follows. In Section \ref{sec:pen_prob} we prove the existence of minimizing functions $\uei{k}$ for $\Iei{k}$ ($k=0,1$) and  show that $\uei{k} \in C^{1,\alpha}(\overline{B})$ for every $\alpha \in (0,1)$. In Section \ref{sec:volcond} we analyze the volume of $\Omega(\uei{0})$ and $\Omega(\uei{1})$, separately. Thereby, Section \ref{sec:nonrew} is concerned with the functional $\Iei{0}$. We will prove Theorem \ref{theo:main} by scaling arguments. In Section \ref{sec:rew} we consider the functional $\Iei{1}$. The rewarding property of the penalization term $\pen{1}$ prohibits an adoption of the scaling argument from Section \ref{sec:nonrew} to prove that $\abs{\Omega(\uei{1})}$ cannot become less than $\omega_0$. Instead, we will use an inequality by M. S. Ashbaugh and R. S. Laugesen (see \cite{AshLau1996}) to establish a lower bound $\alpha_0$ on $\abs{\Omega(\uei{1})}$. This lower bound satisfies $\alpha_0 \in (\frac{1}{2},1)$ and depends on $n,\eps$ and $\omega_0$. 
Thus, choosing $\eps$ sufficiently small, there holds  $\abs{\Omega(\uei{1})} \in [\alpha_0\omega_0,\omega_0]$. Refining the choice of $\eps$ we will obtain the following dichotomy: there either holds
\begin{compactenum}[(i)]\setlength{\itemsep}{0pt}
 \item $\abs{\Omega(\uei{1})} = \omega_0$ \qquad or 
 \item $\abs{\Omega(\uei{1})} < \omega_0$ and if $\Omega(\uei{1})$ is scaled to the volume $\omega_0$, this enlarged domain is neither a  subset of $B$ nor can it be translated into $B$.
\end{compactenum}
In the first case, $\uei{1}$ minimizes the functional $\Iei{0}$ and we may treat $\Iei{0}$ and $\Iei{1}$ as equivalent in the sense that they are minimized by the same functions. 

Under the additional assumption that $\partial\Omega(\uei{1})$ satisfies a doubling condition, we will disprove the occurrence of the latter case. Assuming that there exists a $\sigma>0$ and an $R_0>0$ such that for every $x_0 \in \partial\Omega(\uei{1})$ and every $0<R\leq R_0$ there holds 
\[
  \abs{B_{2R}(x_0)\cap \Omega(\uei{1})} \leq \sigma \abs{B_{R}(x_0)\cap \Omega(\uei{1})},
\]
we will be able to establish a nondegeneracy result for $\uei{1}$. The nondegeneracy will eventually disprove the second case of the above dichotomy. Besides the assumption of the doubling property, the rewarding part of $\pen{1}$ will be crucial for proving the nondegeneracy.

\section{The penalized problems}\label{sec:pen_prob}
We start with proving the existence of solutions for the penalized problems 
\[
   \min\{\Iei{k}(v): v \in H^{2,2}_0(B)\}
\]
for $k=0$ and $k=1$.   
The proof is exactly done as in \cite[Theorem 2.1]{Sto2015} with some obvious changes. 
\begin{theorem}
 For $k \in \{0,1\}$ and every $\eps>0$ there exists an $\uei{k} \in H^{2,2}_0(B)$ such that 
 \[
  \Iei{k}(\uei{k}) = \min\{\Iei{k}(v) : v\in H^{2,2}_0(B)\}.
 \]
\end{theorem}

 Without loss of generality, we assume that $\uei{k}$ is normalized in the sense that
 \[
   \int\limits_B\ue^2 dx =1
 \]
 and we denote 
\[
  \Gei{k} := \Ray(\uei{k},B) = \int\limits_B \abs{\lap \uei{k}}^2dx. 
\]
Note that $\O(\uei{k}) = \{x\in B:\uei{k}(x)\neq 0\}$ is a non-empty set and the absolute continuity of the Lebesgue intregal implies that the $n$-dimensional Lebesgue measure of $\O(\uei{k})$ cannot vanish. 
However, we do not know if $\O(\uei{k})$ is an open set as long as we do not possess any information about regularity properties of $\uei{k}$.
Our next step is to establish the $C^{1,\alpha}$ regularity of $\ue$. For this purpose, we follow an idea of Q. Han and F. Lin \cite{HanLin}, which is based on Morrey's Dirichlet Growth Theorem (see Theorem \ref{morrey}). 

In the sequel, we will apply the following version of Morrey's Dirichlet Growth Theorem. For the proof we refer to \cite{morrey}.
\begin{theorem}\label{morrey}
 Suppose $\phi \in H^{1,p}_0(B)$, $1 \leq p\leq n$, $0< \alpha\leq1$ and  suppose there exists a constant $M>0$ such that 
 \[
  \int\limits_{B_r(x_0)\cap B}\abs{\grad \phi}^pdx \leq M \, r^{n-p+\alpha p}
 \]
for every $B_r(x_0)$ with $x_0 \in \overline{B}$. Then $\phi \in C^{0,\alpha}(\overline{B})$.
\end{theorem}

We need to verify the assumptions of Theorem \ref{morrey} for the second order derivatives of $\ue$. 
Let us fix an $R_{\eps,k}$ with $0 < R_{\eps,k} <1$ such that the volume of a ball with radius $R_{\eps,k}$ is smaller than $\abs{\O(\ue)}$. 
Now suppose $x_0 \in \overline{B}$ and choose  $0 <r \leq R \leq R_{\eps,k}$. We define $\hat{v}_k \in H^{2,2}_0(B)$ by 
\[
 \hat{v}_k = \begin{cases}
            \ue & \mbox{ in } B\setminus B_R(x_0) \\
            v_k & \mbox{ in } B_R(x_0) \cap B
           \end{cases},
\]
where $v_k-\ue \in H^{2,2}_0(B_R(x_0)\cap B)$ and $\bilap v_k =0$ in $B_R(x_0)\cap B$.
 Then the estimate
\begin{equation}\label{eq:cont_1}
 \int\limits_{B_r(x_0)\cap B}\abs{D^2 \ue}^2dx \leq 2 \int\limits_{B_R(x_0)\cap B}\abs{D^2(\ue-v_k)}^2dx + 2 \int\limits_{B_r(x_0)\cap B}\abs{D^2v_k}^2dx
\end{equation}
is obvious. The next lemma helps to estimate the last term in the above inequality.
\begin{lemma}\label{reg_1_la_1}
Using the above notation, there exists a constant $C>0$ such that for each  $r \leq R$ the following estimate holds
 \[
   \int\limits_{B_r(x_0)\cap B}\abs{D^2 v_k}^2dx \leq C\,\left(\frac{r}{R}\right)^n \,\int\limits_{B_R(x_0)\cap B}\abs{D^2\ue}^2dx\,.
 \]
Thereby, the constant $C$ is independent of $r,\,R$ and $x_0$\,.
\end{lemma}
The proof of this lemma is analogue to \cite[Lemma 2.1]{Sto2015} with some obvious changes. 

It remains to estimate the first integral on the right hand side of \eqref{eq:cont_1}. For this purpose, we need the following  lemmata. The first one is a direct consequence of the minimality of $\ue$ for $\Iei{k}$. 

\begin{lemma}\label{la:tech_morrey2}
  There exists a constant $C= C(n,\Ge,\abs{\O(\ue)})>0$ such that for each $x_0 \in \overline{B}$ and each $0<R\leq R_{\eps,k}$ there holds
 \[
  \int\limits_{B_R(x_0)\cap B}\abs{D^2(\ue-v_k)}^2 dx \leq C(n,\Ge,\abs{\O(\ue)})\left (R^n + \int\limits_{B_R(x_0)} \ue^2dx \right).
 \]
\end{lemma}
\begin{proof} 
Note that $\O(\ue) \not\subset B_R(x_0)$ due to the definition of $R_{\eps,k}$. However, we have to distinguish two different cases relating to the volume of $\O(\hat{v}_k)$. First, let us consider that there holds $\abs{\O(\hat{v}_k)} \leq \abs{\O(\ue)}$. 
Then the minimality of $\ue$ for $\Iei{k}$ implies
 \begin{align*}
  \Ge \int\limits_B \hat{v}_k^2dx &\leq \int\limits_B \abs{\lap \hat{v}_k}^2dx 
  \aq \int\limits_{B_R(x_0)\cap B}\abs{\lap\ue}^2-\abs{\lap{v_k}}^ 2dx &\leq \Ge \int\limits_{B_R(x_0)\cap B}\ue^2-v_k^2dx
 \end{align*}
Thus,
\[
 \int\limits_{B_R(x_0)\cap B}\abs{D^2(\ue-v_k)}^2dx \leq \Ge \int\limits_{B_R(x_0)}\ue^2dx.
\]
Now let us assume that there holds $\abs{\O(\ue)} < \abs{\O(\hat{v}_k)}$.
In this case, the penalization term inhibits a reasonable comparison of $\Iei{k}(\ue)$ and $\Iei{k}(\hat{v})$. We circumvent this problem by scaling $\hat{v}_k$. Let $R_B$ be the radius of the reference domain $B$. Then we choose a ball $B^\ast$, which is concentric to $B$. The radius of $B^\ast$ is 
\[
 R^\ast := \frac{R_B}{\mu}, \quad \mbox{ where } \mu := \left(\frac{\abs{\O(\ue)}}{\abs{\O(\hat{v}_k)}}\right)^\frac{1}{n}.
\]
Note that $\mu>1$. Now we set $w_k(x) := \hat{v}_k(\mu x)$ for $x\in B^\ast$. Note that $w_k \in H^{2,2}_0(B^\ast)\subset H^{2,2}_0(B)$ and that $\abs{\O(w_k)}=\abs{\O(\ue)}$. Additionally, we find
\begin{equation}\label{eq:mu}
1 < \mu \leq \left( 1+\frac{\abs{B_R(x_0)}}{\abs{\O(\ue)}} \right)^\frac{1}{n} \quad \mbox{ and } \quad 
  1-\mu^{-4} \leq C(n,\abs{\O(\ue)})R^n.
\end{equation}
The minimality of $\ue$ for $\Iei{k}$ in $H^{2,2}_0(B)$ now implies
\[
 \Ge \int\limits_{B}w_k^2dx \leq \int\limits_{B}\abs{\lap w_k}^2dx 
 \quad \aq \quad  \Ge \mu^{-4}\int\limits_{B}\hat{v}_k^2dx \leq \int\limits_{B}\abs{\lap\hat{v}_k}^2dx .
\]
Rearranging terms we obtain the local inequality
\begin{align*}
 \int\limits_{B_R(x_0)\cap B}\abs{\lap \ue}^2-\abs{\lap v_k}^2dx \leq \Ge\left(1-\mu^{-4}\right)+\frac{\Ge}{\mu^4}\int\limits_{B_R(x_0)\cap B}\ue^2-v_k^2dx
\end{align*}
Using the second estimate in \eqref{eq:mu}, we find that
\begin{align*}\label{eq:est_n>4}
 \int\limits_{B_R(x_0)\cap B} \abs{D^2(\ue-v_k)}^2dx\leq C(n,\Ge,\abs{\O(\ue)}) \left(R^n + \int\limits_{B_R(x_0)}\ue^2dx\right).
\end{align*}
\end{proof}

The following technical lemma is cited from \cite[Chapter III, Lemma 2.1]{giaq_mult_int}.
\begin{lemma}\label{la:tech_morrey}
 Let $\Phi$ be a nonnegative and nondecreasing function on $[0,R]$. Suppose that there exist positive constants $\gamma, \alpha, \kappa, \beta$, $\beta<\alpha$, such that for all $0\leq r\leq R \leq R_0$ 
\[
  \Phi(r) \leq \gamma\left[ \left(\frac{r}{R}\right)^\alpha + \delta \right]\Phi(R)+\kappa R^\beta.
\]
Then there exist positive constants $\delta_0 =\delta_0(\gamma,\alpha,\beta)$ and $C=C(\gamma, \alpha,\beta)$ such that if $\delta < \delta_0$, for all $0\leq r\leq R\leq R_0$ we have
\[
  \Phi(r) \leq C \left(\frac{r}{R}\right)^\beta \left[\Phi(R) + \kappa R^\beta\right]  .
\] 
\end{lemma}

The next lemma is the essential tool for proving the $C^{1,\alpha}$ regularity of $\ue$. It is based on ideas of \cite[Chapter 3]{HanLin}. 

\begin{lemma}\label{la:bootstrap1}
 Suppose that for each $0\leq r\leq  R_{\eps,k}$ there holds
 \[
  \int\limits_{B_r(x_0)} \abs{D^2 \ue}^2dx \leq M\, r^\mu, \quad \mu \in [0,n). 
 \]
Then there exists a constant $C(n,\abs{\O(\ue)})>0$ such that for each $0\leq r \leq R_{\eps,k}$ 
\[
 \int\limits_{B_r(x_0)}\ue^2dx \leq C(n,\Ge,\abs{\O(\ue)})(1+M)\,r^\sigma
\]
where $\sigma = \mu +4$ if $\mu < n-4$ and $\sigma$ is arbitrary in $(0,n)$ if $n-4 \leq \mu < n$.
\end{lemma}
\begin{proof}
 Let $0\leq r\leq s \leq R_{\eps,k}$. For a function $w \in H^{2,2}(B)$ we set 
 \[
  (w)_{r,x_0} := \fint\limits_{B_r(x_0)}w\,dx.
 \]
We split the proof in two parts. First, we estimate $\norm{\grad \ue}_{L^2}^2$.  
 Young's inequality implies
\begin{align*}
 \int\limits_{B_r(x_0)}\abs{\grad \ue}^2dx &\leq 2 \sum_{i=1}^n \left( \int\limits_{B_r(x_0)}(\partial_i\ue)^2_{s,x_0}dx + \int\limits_{B_r(x_0)}\abs{\partial_i\ue - (\partial_i\ue)_{s,x_0}}^2dx \right) \\
 &\leq 2 \sum_{i=1}^n \left(\abs{B_r} \left(\int\limits_{B_s(x_0)}\partial_i\ue\,dx \right)^2 + \int\limits_{B_s(x_0)}\abs{\partial_i\ue - (\partial_i\ue)_{s,x_0}}^2dx \right).
\end{align*}
Applying H\"older's and Poincaré's inequality, we find that
\begin{align*}
 \int\limits_{B_r(x_0)}\abs{\grad \ue}^2dx &\leq C \left[ \left(\frac{r}{s}\right)^n \int\limits_{B_s(x_0)}\abs{\grad \ue}^2dx + s^2\int\limits_{B_s(x_0)}\abs{D^2\ue}^2dx  \right].
\end{align*}
By assumption, we can proceed to
\begin{align*}
 \int\limits_{B_r(x_0)}\abs{\grad \ue}^2dx &\leq C \left[\left(\frac{r}{s}\right)^n \int\limits_{B_s(x_0)}\abs{\grad \ue}^2dx +M\, s^{\mu+2} \right].
\end{align*}
Now Lemma \ref{la:tech_morrey} implies that for each $0\leq r\leq s\leq R_{\eps,k}$ there holds
\begin{align*}
 \int\limits_{B_r(x_0)}\abs{\grad \ue}^2dx &\leq C  \left(\frac{r}{s}\right)^\theta \left[\int\limits_{B_{s}(x_0)}\abs{\grad \ue}^2dx + M\,s^\theta  \right],
\end{align*}
where $\theta = \mu +2$ if $\mu< n-2$ and $\theta$ is arbitrary in $(0,n)$ if $n-2 \leq \mu < n$. Since $\norm{\ue}_{H^{2,2}_0(B)}=\Ge$, we obtain
\[ 
\int\limits_{B_r(x_0)}\abs{\grad \ue}^2dx    \leq C  \left(\frac{r}{s}\right)^\theta \left[\Ge + M\,s^\theta  \right].
\]
Choosing $s=R_{\eps,k}$, we deduce
\begin{equation}\label{eq:aux_1}
 \int\limits_{B_r(x_0)}\abs{\grad \ue}^2dx \leq C(n,\abs{\O(\ue)})(\Ge+M)r^\theta.
\end{equation}
We estimate $\norm{\ue}^2_{L^2}$ in the same way. Let $0<r \leq s \leq R_{\eps,k}$, then
\begin{align*}
 \int\limits_{B_r(x_0)}\ue^2dx &\leq 2 \int\limits_{B_r(x_0)}(\ue)^2_{s,x_0}dx + 2\int\limits_{B_s(x_0)}\abs{\ue-(\ue)_{s,x_0}}^2dx \\ 
 &\leq C\left(\frac{r}{s}\right)^n\int\limits_{B_s(x_0)}\ue^2dx + C\,s^2\,\int\limits_{B_s(x_0)}\abs{\grad \ue}^2dx.
\end{align*}
Hence, applying \eqref{eq:aux_1} there holds 
\[
 \int\limits_{B_r(x_0)}\ue^2dx \leq C\left(\frac{r}{s}\right)^n\int\limits_{B_s(x_0)}\ue^2dx + C(n,\abs{\O(\ue)})(\Ge+M)s^{\theta+2}
\]
and by Lemma \ref{la:tech_morrey} we obtain  
\begin{align*}
  \int\limits_{B_r(x_0)}\ue^2dx &\leq C\left(\frac{r}{s}\right)^\sigma\left( \int\limits_{B_{s}(x_0)}\ue^2dx + C(n,\abs{\O(\ue)})(\Ge+M)s^\sigma\right), 
\end{align*}
where $\sigma = \theta + 2$ if $\theta < n-2$ and $\sigma$ arbitrary in $(0,n)$ if $n-2 \leq \theta < n$. 
Thus, choosing $s=R_{\eps,k}$, there holds
\[
  \int\limits_{B_r(x_0)}\ue^2dx \leq C(n,\abs{\O(\ue)})(1+\Ge+M)r^\sigma \leq C(n,\Ge,\abs{\O(\ue)})(1+M)r^\sigma.
\]
for each $0\leq r\leq R_{\eps,k}$. This proves the claim.
\end{proof}

We are now ready  to prove the $C^{1,\alpha}$ regularity of $\ue$. 

\begin{theorem}\label{theo:reg1}
  For each $\eps>0$ and $k\in\{0,1\}$ the minimizer $\ue$ suffices $\ue \in C^{1,\alpha}(\overline{B})$ for every $\alpha \in (0,1)$.  
\end{theorem}
\begin{proof}
Our aim is  to show that for each $x_0 \in \overline{B}$ and every $0<r\leq R_{\eps,k}$ there holds 
 \begin{equation}\label{eq:boot_0}
  \int\limits_{B_r(x_0)}\abs{D^2\ue}^2 dx \leq C(n,\Ge,\abs{\O(\ue)})\,r^{n-2+2\alpha}
 \end{equation}
for each $\alpha \in (0,1)$. Then Theorem \ref{morrey} finishes the proof if we choose $\phi = \partial_i\ue$ for $i\in\{1,\ldots,n\}$. 
We prove \eqref{eq:boot_0} by using a bootstrap argument. 
 Let $x_0 \in \overline{B}$ and $0\leq r\leq R\leq R_{\eps,k}$. Due to estimate \eqref{eq:cont_1}, Lemma \ref{reg_1_la_1} and Lemma \ref{la:tech_morrey2} there holds
 \begin{equation}\label{eq:boot0}\begin{split}
   &\int\limits_{B_r(x_0)\cap B}\abs{D^2\ue}^2dx \\
  \leq &\quad C \left(\frac{r}{R}\right)^n\int\limits_{B_R(x_0)}\abs{D^2\ue}^2dx + C(n,\Ge,\abs{\O(\ue)})\left(R^n + \int\limits_{B_R(x_0)}\ue^2dx\right).\end{split}
 \end{equation}
Thus, for every $0\leq r\leq R_{\eps,k}$ there holds
\begin{equation}\label{eq:boot1}
  \int\limits_{B_r(x_0)\cap B}\abs{D^2\ue}^2dx \leq C(n,\Ge,\abs{\O(\ue)})\,r^0.
\end{equation}
Consequently, for every $0\leq r \leq R_{\eps,k}$, Lemma \ref{la:bootstrap1} implies 
\begin{equation*}
 \int\limits_{B_r(x_0)}\ue^2dx \leq C(n,\Ge,\abs{\O(\ue)})\,r^{\sigma_0},
\end{equation*}
where $\sigma_0 = 4$ if $n\geq 5$ and $\sigma_0$ arbitrary in $(0,n)$ if $n=4$. We insert this estimate in \eqref{eq:boot0}. This yields
\begin{align*}
 \int\limits_{B_r(x_0)\cap B}\abs{D^2\ue}^2dx &\leq C \left(\frac{r}{R}\right)^n\int\limits_{B_R(x_0)}\abs{D^2\ue}^2dx + C(n,\Ge,\abs{\O(\ue)})\left(R^n + R^{\sigma_0}\right) \\
 &\leq C \left(\frac{r}{R}\right)^n\int\limits_{B_R(x_0)\cap B}\abs{D^2\ue}^2dx + C(n,\Ge,\abs{\O(\ue)})R^{\sigma_0}
\end{align*}
for every $0\leq r \leq R\leq R_{\eps,k}$. Applying Lemma \ref{la:tech_morrey}, we obtain
\begin{align*}
 \int\limits_{B_r(x_0)\cap B}\abs{D^2\ue}^2dx &\leq C \left(\frac{r}{R}\right)^{\sigma_0}\left(\int\limits_{B_R(x_0)\cap B}\abs{D^2\ue}^2dx + C(n,\Ge,\abs{\O(\ue)})R^{\sigma_0}\right) 
\end{align*}
for $0\leq r\leq R\leq R_{\eps,k}$. Then choosing $R =R_{\eps,k}$ gives us
\begin{equation}\label{eq:boot2}
  \int\limits_{B_r(x_0)\cap B}\abs{D^2\ue}^2dx \leq  C(n,\Ge,\abs{\O(\ue)})\,r^{\sigma_0}
\end{equation}
for every $0\leq r\leq R_{\eps,k}$. Note that \eqref{eq:boot2} is an improvement of estimate \eqref{eq:boot1}. 

For $n=4$, estimate \eqref{eq:boot2} proves the claim. For $n\geq 5$, we repeat the argumentation. Based on \eqref{eq:boot2}, Lemma \ref{la:bootstrap1}  implies that 
\[
 \int\limits_{B_r(x_0)}\ue^2dx \leq C(n,\Ge,\abs{\O(\ue)})\,r^{\sigma_1},
\]
where $\sigma_1 = 8 $ if $n\geq 8$ and $\sigma_1$ is arbitrary in $(0,n)$ if $5\leq n\leq 7$. Again, we insert this estimate in \eqref{eq:boot0} and deduce an improvement of \eqref{eq:boot2}. Repeating this process proves the claim after finite many steps.
\end{proof}

Due to Theorem \ref{theo:reg1}, the set $\O(\ue)$ is an open set and classical variational arguments show that 
\[
  \bilap \ue - \Ge \ue =0\; \mbox{ in  } \O(\ue) .
\] 
Moreover, the $C^{1,\alpha}$ regularity of $\ue$ allows us to split $\partial\O(\ue)$ in the following two parts
\[
 \Sigma_{\eps,k}^0 := \{x \in \partial\O(\ue): \abs{\grad\ue(x)}=0\} \; \mbox{ and } \; \Sigma_{\eps,k}^1 := \{x \in \partial\O(\ue): \abs{\grad\ue(x)}>0\}.
\]
Then $\Sigma^1_{\eps,k}$ is part of a nodal line of $\ue$ since $\ue \in H^{2,2}_0(B)$ and, consequently, the $n$-dimensional Lebesgue measure of $\Sigma^1_{\eps,k}$ vanishes for every $\eps>0$. 
We define
\[
 \Omega(\ue) := \O(\ue) \cup \Sigma^1_{\eps,k}
\]
and call $\Sigma^0_{\eps,k} = \partial\Omega(\ue)$ the free boundary. 

\begin{remark}
  Note that $\Omega(\ue)$ is an open set in $\R^n$ and $\abs{\O(\ue)}=\abs{\Omega(\ue)}$. Moreover, the minimizer $\ue$ solves
  \[
    \begin{cases}
      \bilap \ue - \Ge \ue &=0\; \mbox{ in } \Omega(\ue) \\
   \ue = \abs{\grad\ue} &= 0 \; \mbox{ on }\partial\Omega(\ue).
    \end{cases}
  \]
\end{remark}

The following lemma shows that, considering the functional $\Iei{1}$, the set $\Omega(\uei{1})$ is connected and, hence, $\Omega(\uei{1})$ is a domain. This result is a direct consequence of the strict monotonicity of the penalization term $\pen{1}$.
\begin{lemma}\label{la:connected_1}
  For every $\eps>0$ the set $\Omega(\uei{1})$ is connected. 
\end{lemma}
\begin{proof}
  We prove the claim by contradiction. Let us assume that $\Omega(\uei{1})$ consists of two connected components $\Omega_1$ and $\Omega_2$ with $\abs{\Omega_k}\neq 0$ for $k=1,2$. Moreover, we define for $k=1,2$
  \[
     u_k :=  \begin{cases}
       \uei{1}, &\mbox{ in } \Omega_k \\
       0, &\mbox{ otherwise }. 
     \end{cases}
  \]
The minimality of $\uei{1}$ for $\Iei{1}$ implies 
\[
   \Iei{1}(\uei{1}) = \Gei{1} + \pen{1}(\abs{\Omega(\uei{1})}) \leq \Iei{1}(u_1) = \Ray(u_1) + \pen{1}(\abs{\Omega_1}).
\]
The definition of $\Ray(u_1)$ and the normalization of $\ue$ now give 
\begin{align*}
  \left( \int\limits_{\Omega_1}\abs{\lap u_1}^2dx + \int\limits_{\Omega_2}\abs{\lap u_2}^2dx\right)&\left(1 - \int\limits_{\Omega_2}u_2^2dx\right) \\ 
  \leq \int\limits_{\Omega_1}&\abs{\lap u_1}^2dx + (\pen{1}(\abs{\Omega_1}) - \pen{1}(\abs{\Omega(\uei{1})}))\int\limits_{\Omega_1}u_1^2 dx
\end{align*}
Since the penalization term $\pen{1}$ is strictly increasing we deduce
\begin{equation}\label{eq:ray_u2}
  \int\limits_{\Omega_2}\abs{\lap u_2}^2dx < \Gei{1}\int\limits_{\Omega_2}u_2^2dx \; \aq \Ray(u_2)<\Gei{1}\,
\end{equation}
We again make use of the minimality of $\uei{1}$ for $\Iei{1}$ and obtain
\[
   \Iei{1}(\uei{1}) = \Gei{1} + \pen{1}(\abs{\Omega(\uei{1})}) \leq \Ray(u_2) + \pen{1}(\abs{\Omega_2}). 
\]
Then \eqref{eq:ray_u2} yields
\[
  \pen{1}(\abs{\Omega(\uei{1})}) = \pen{1}(\abs{\Omega_1}+\abs{\Omega_2}) < \pen{1}(\Omega_2).
\]
Since $\pen{1}$ is strictly increasing, we derive a contradiction and the claim is proven. 
\end{proof}

Note that for the proof of Lemma \ref{la:connected_1} the actual value of $\abs{\Omega(\uei{1})}$ is irrelevant since $\pen{1}$ is stricly increasing. Considering the functional $\Iei{0}$ with the only nondecreasing penalization term $\pen{0}$, we have to ensure that $\abs{\Omega(\uei{0})}\geq \omega_0$ before we are able to copy the approach of Lemma \ref{la:connected_1} and obtain that $\Omega(\uei{0})$ is connected. This will be done in Section \ref{sec:nonrew}. 

The next corollary collects direct consequences of Lemma \ref{la:connected_1}.
\begin{corollary}\label{cor:optdomain_1}
  For every minimizer $\uei{1}$ of $\Iei{1}$ the domain $\Omega(\uei{1})$ satisfies $\Gei{1}=\Gamma(\Omega(\uei{1}))$ and $\Omega(\uei{1})$ is an optimal domain for minimizing the fundamental tone among all open subset of $B$ with the same measure as $\Omega(\uei{1})$. 
\end{corollary}
\begin{proof}
 Let $\uei{1} \in H^{2,2}_0(B)$ minimize $\Iei{1}$. Then there holds 
 \begin{align*}
    \Gei{1} + \pen{1}(\abs{\Omega(\uei{1})}) &= \min\{\Iei{1}(v) : v \in H^{2,2}_0(B)\} \\
    & \leq \min\{ \Ray(v) + \pen{1}(\abs{\O(v)}) : v \in H^{2,2}_0(B), \abs{\O(v)} = \abs{\Omega(\uei{1})}\} \\
    & \leq \min\{ \Ray(v)  : v \in H^{2,2}_0(B)\}+ \pen{1}(\abs{\Omega(\uei{1})}) \\
    & = \G(\Omega(\uei{1})) +  \pen{1}(\abs{\Omega(\uei{1})}) \\
    & \leq \Gei{1} + \pen{1}(\abs{\Omega(\uei{1})}).
 \end{align*}
 Thus, $\Gei{1} = \G(\Omega(\uei{1}))$ and $\uei{1}$ is a clamped eigenfunction on $\Omega(\uei{1})$. 
\end{proof}

\section{The volume condition}\label{sec:volcond}

\subsection{The non-rewarding penalization term}\label{sec:nonrew}
In this section, we focus on the functional $\Iei{0}$.

\begin{theorem}\label{theo:vol_0_1}
  For every $\eps>0$ there holds $\abs{\Omega(\uei{0})}\geq \omega_0$ for every minimizer $\uei{0}$ of $\Iei{0}$.
\end{theorem}
\begin{proof}
 We prove the claim by contradiction. Therefore, we assume that for an $\eps>0$ there exists a minimizer $\uei{0}$ of $\Iei{0}$ such that $\abs{\Omega(\uei{0})}<\omega_0$. Let us choose an $x_0 \in \partial\Omega(\uei{0})\setminus \partial B$ and a radius $r>0$ such that 
 \begin{equation}\label{eq:vol_0_1_0}
  \abs{B_r(x_0)\cap \{x \in B: \uei{0}(x)\neq 0\}}> 0
 \end{equation}
 and $\abs{\Omega(\uei{0})\cup B_r(x_0)} \leq \omega_0$. Note that such an $x_0 \in \partial\Omega(\uei{0})$ exists since we assume $\abs{\Omega(\uei{0})}<\omega_0$. In addition, let $v \in H^{2,2}(B_r(x_0))$ with $v-\uei{0} \in H^{2,2}_0(B_r(x_0))$ satisfying
 \[
  \bilap v - \Gei{0}v =0 \mbox{ in } B_r(x_0). 
 \]
 We set 
 \[
   \hat{v} := 
   \begin{cases}
       \uei{0}, &\mbox{ in } B\setminus B_r(x_0) \\
       v, & \mbox{ in } B_r(x_0)
   \end{cases}
 \]
 and compare the $\Iei{0}$-energies of $\uei{0}$ and $\hat{v}$. This leads to the following local inequality
 \begin{equation}\label{eq:vol_0_1_1}
    \int\limits_{B_r(x_0)}\abs{\lap \uei{0}}^2dx -\abs{\lap v}^2dx \leq \Gei{0} \int\limits_{B_r(x_0)}\uei{0}^2-v^2dx. 
 \end{equation}
 Applying integration by parts and the definition of $v$, we obtain
 \[
   \int\limits_{B_r(x_0)}\abs{\lap \uei{0}}^2dx -\abs{\lap v}^2dx = \int\limits_{B_r(x_0)}\abs{\lap(\uei{0}-v)}^2dx + 2\,\Gei{0}\int\limits_{B_r(x_0)}v(\uei{0}-v)dx.
 \]
 Thus, \eqref{eq:vol_0_1_1} becomes 
 \[
   \int\limits_{B_r(x_0)}\abs{\lap(\uei{0}-v)}^2dx \leq \Gei{0}\int\limits_{B_r(x_0)}\abs{\uei{0}-v}^2dx
 \]
 and applying Poincaré's inequality twice yields 
 \begin{equation}\label{eq:vol_0_1_2}
  \int\limits_{B_r(x_0)}\abs{\lap(\uei{0}-v)}^2dx \leq r^4\,\Gei{0}\,\int\limits_{B_r(x_0)}\abs{\lap(\uei{0}-v)}^2dx.
 \end{equation}
 Provided that the integral in \eqref{eq:vol_0_1_2} does not vanish, \eqref{eq:vol_0_1_2} is contradictory for sufficently small $r$. 

If the integral in \eqref{eq:vol_0_1_2} vanishes, there holds $\uei{0}\equiv v$ in $B_r(x_0)$. Consequently, $\uei{0}$ is analytic in $B_\frac{r}{2}(x_0)$ since $v$ is there analytic as a solution of an ellipitc equation. However, then $\uei{0}$ vanishes in $B_\frac{r}{2}(x_0)$ because of \eqref{eq:vol_0_1_0}. This is contradictory since $x_0\in\partial\Omega(\uei{0})$. This proves the claim.
\end{proof}

As a consequence of Theorem \ref{theo:vol_0_1} we can adopt the proof of Lemma \ref{la:connected_1} to show that $\Omega(\uei{0})$ is connected. 

\begin{lemma}\label{la:connected_0}
  For every $\eps>0$ and every minimizer $\uei{0}$ of $\Iei{0}$ the set $\Omega(\uei{0})$ is connected. 
\end{lemma}
\begin{proof}
 We follow the idea of the proof of Lemma \ref{la:connected_1} and assume that $\Omega(\uei{0})=\Omega_1 \dot{\cup} \Omega_2$ with $\abs{\Omega_k}>0$ for $k = 1,2$.  Let us first consider that $\abs{\Omega(\uei{0})}>\omega_0$. Then there holds
 \begin{equation}\label{eq:conn_1}
   \pen{0}(\abs{\Omega(\uei{0})}) > \pen{0}(\Omega_k) 
 \end{equation}
 for $k= 1,2$. Hence, arguing as in the proof of Lemma \ref{la:connected_1} and replacing the strict monotonicity of $\pen{1}$ by \eqref{eq:conn_1}, we arrive at a contradiction. 
 
 Now let us assume that $\abs{\Omega(\uei{0})}=\omega_0$. Then 
  \begin{equation}\label{eq:conn_2}
   \pen{0}(\abs{\Omega(\uei{0})}) = \pen{0}(\Omega_k) =0
 \end{equation}
 and arguing as in the proof of Lemma \ref{la:connected_1} we obtain 
 \[
   \Iei{0}(u_2)=\Ray(u_2) \leq \Gei{0} = \Iei{0}(\uei{0}).
 \]
 This implies that $u_2$ minimizes the functional $\Iei{0}$. Since $\Omega(u_2)=\Omega_2$ and $\abs{\Omega_2}<\omega_0$ this is contradictory to Theorem \ref{theo:vol_0_1}. Thus, the claim is proven.
\end{proof}

As a consequence of Lemma \ref{la:connected_0}, we get the analog to Corollary \ref{cor:optdomain_1}. 
\begin{corollary}\label{cor:optdomain_0}
For every minimizer $\uei{0}$ of $\Iei{0}$ the set $\Omega(\uei{0})$ is a domain and $\Gei{0} = \G(\Omega(\uei{0}))$. In addition, $\Omega(\uei{0})$ minimizes the fundamental tone among all open subsets of $B$ with the same measure as $\Omega(\uei{0})$.
\end{corollary}   

The following remark will be helpful to show that for an appropriate choice of $\eps$ the $n$-dimensional Lebesgue measure of $\Omega(\uei{0})$ cannot become larger than $\omega_0$. 

\begin{remark}\label{rem:EF_omega_0}
 Note that the reference domain $B$ compactly contains a ball $B_R(x_0)$ with $\abs{B_R}=\omega_0$. Let $\phi\in H^{2,2}_0(B_R(x_0))$ denote the clamped eigenfunction on $B_R(x_0)$, i.e.
 \[
  \G(B_R(x_0)) = \min\{\Ray(v): v\in H^{2,2}_0(B_R(x_0))\} = \Ray(\phi).
 \]
 Consequently, for every $\eps>0$ and $k\in\{0,1\}$ there holds
\[
  \Iei{k}(\phi) = \Ray(\phi) = \G(B_R(x_0)) = \left(\frac{\omega_n}{\omega_0}\right)^\frac{4}{n}\G(B_1),
\]
where $B_1$ denotes the unit ball in $\R^n$.
\end{remark}

\begin{theorem}\label{theo:vol_0_2}
Let $\uei{0}$ be a minimizer for $\Iei{0}$. Then there exists a number $\eps_1 = \eps_1(n,\omega_0)$ such that 
\[
 \abs{\Omega(\uei{0})} \leq \omega_0.
\]
\end{theorem}
\begin{proof}
We claim that the statement of the theorem holds true for 
\begin{equation}\label{eq:eps_1}
  \eps_1 = \eps_1(n,\omega_0) = \left( \frac{\omega_0}{\omega_n}\right)^\frac{4}{n}\,\frac{\omega_0}{\G(B_1)}. 
\end{equation}
We prove by contradiction. Thus, we choose $\eps\leq\eps_1$ and denote by $\uei{0}$ a minimizer of $\Iei{0}$. Assume that there exists a number $\alpha>1$
 such that 
\begin{equation}\label{eq:assum_vol_larger}
  \abs{\Omega(\uei{0})} = \alpha \omega_0.
\end{equation}
Our aim is to contradict \eqref{eq:assum_vol_larger}. 
Since $\Omega(\uei{0}) \subset B$ and $\alpha >1$, the scaled domain $\Omega' := \alpha^{-\frac{1}{n}}\Omega(\uei{0})$ is also contained in $B$ and satisfies $\abs{\Omega'}=\omega_0$. Let $\psi \in H^{2,2}_0(B)$ denote an eigenfunction for the bi-Laplacian on $\Omega'$. Then the minimality of $\uei{0}$ for $\Iei{0}$ implies
\begin{equation}\label{eq:vol_larger_0}
   \G(\Omega(\uei{0})) + \frac{\omega_0}{\eps}(\alpha -1) = \Iei{0}(\uei{0}) \leq \Iei{0}(\psi) = \G(\Omega').
\end{equation}
By scaling we have 
\[
  \G(\Omega') = \G(\alpha^{-\frac{1}{n}}\Omega(\uei{0})) = \alpha^\frac{4}{n}\G(\Omega(\uei{0})). 
\] 
From \eqref{eq:vol_larger_0} we then get
\begin{equation}\label{eq:vol_larger_1}
   \frac{\omega_0}{\eps}(\alpha -1) \leq \left(\alpha^\frac{4}{n}-1\right)\G(\Omega(\uei{0})). 
\end{equation}
Now let $\phi \in H^{2,2}_0(B)$ be as in Remark \ref{rem:EF_omega_0}. Then, due to the assumption \eqref{eq:assum_vol_larger}, there holds
\[
   \G(\Omega(\uei{0})) < \Iei{0}(\uei{0}) \leq \Iei{0}(\phi) = \left(\frac{\omega_n}{\omega_0}\right)^\frac{4}{n}\G(B_1)
\]
and estimate \eqref{eq:vol_larger_1} becomes
\[
   \left(\frac{\omega_0}{\omega_n}\right)^\frac{4}{n}\,\frac{\omega_0}{\G(B_1)}\,\frac{\alpha -1}{\alpha^\frac{4}{n}-1}< \eps. 
\]
Hence, 
\[
  \frac{\alpha-1}{\alpha^\frac{4}{n}-1} < \frac{\eps}{\omega_0}\,\left(\frac{\omega_n}{\omega_0}\right)^\frac{4}{n}\G(B_1). 
\]
With \eqref{eq:eps_1} this implies
\[
1 \leq  \frac{\alpha-1}{\alpha^\frac{4}{n}-1} < \eps\,\eps_1^{-1}
\]
 since $\alpha >1$ and $n\geq 4$. Thus, for any $\eps\leq\eps_1$ we get a contradiction and the assumption \eqref{eq:assum_vol_larger} cannot hold true. Together with Theorem \ref{theo:vol_0_1} the claim is proven.
\end{proof}

Finally, the proof of the main theorem, Theorem \ref{theo:main}, is a direct consequence of the previous results in this section. 

\begin{proof}[Proof of Theorem \ref{theo:main}]
   Let us choose $\eps\leq\eps_1$, where $\eps_1$ is given in Theorem \ref{theo:vol_0_2} and let $\uei{0}$ be a minimizer of $\Iei{0}$. Then Theorem \ref{theo:vol_0_2} implies that $ \abs{\Omega( \uei{0})}=\omega_0$ and, due to Corollary \ref{cor:optdomain_0}, there holds 
   \[
      \Iei{0}(\uei{0})=\G(\Omega(\uei{0})). 
   \]
   Now choose an open set $D \subset B$ with $\abs{D} \leq \omega_0$ and denote by $u_D \in H^{2,2}_0(D)$ the clamped eigenfunction on $D$. Then the minimality of $\uei{0}$ for $\Iei{0}$ implies 
   \[
     \G(\Omega(\uei{0})) = \Iei{0}(\uei{0}) \leq \Iei{0}(u_D) = \G(D). 
   \]
   Hence,
   \begin{equation}\label{eq:opt1}
      \G(\Omega(\uei{0})) = \min\{\G(D) : D\subset B, D \mbox{ open, } \abs{D} \leq \omega_0\}.
   \end{equation}
   In addition, $\Omega(\uei{0})$ is connected (see Lemma \ref{la:connected_0}). This proves Theorem \ref{theo:main}.
\end{proof}

Since the existence of an optimal domain for minimizing $\G$  among all open subsets of $B$ of given volume is now proven, the next reasonable step would be a qualitative analysis of the free boundary $\partial\Omega(\uei{0})$. 
Following \cite{AltCaf81}, our next aims would by establishing a nondegeneracy result for $\uei{0}$. Considering second order problems (e.g. \cite{AguAltCaf1986,AltCaf81,BaWa09}), these nondegeneracy results are achieved by applying comparision principles, which are not available for fourth order operators in general. One possible way out of this difficulty is to replace the penalization term $\pen{0}$ by the rewarding penalization term $\pen{1}$. This will be discussed in the next section. 

\subsection{The rewarding penalization term} \label{sec:rew}

In this section, we consider the functional $\Iei{1}$. Analog to Theorem \ref{theo:vol_0_2} we will find that $\abs{\Omega(\uei{1})}$ cannot become larger than $\omega_0$ provided that $\eps\leq\eps_1$. 

It remains to exclude that $\Omega(\uei{1})<\omega_0$. This will be more involved since adopting a scaling argument like the one we used in the proof of Theorem \ref{theo:vol_0_2} collapses if we cannot guarantee that the scaled version of $\Omega(\uei{1})$ is still contained in the reference domain $B$. 

Choosing the parameter $\eps$ sufficiently small, we will see that one of the following two situations occurs: either $\abs{\Omega(\uei{1})}=\omega_0$ or $\abs{\Omega(\uei{1})}<\omega_0$ and the rescaled domain $\Omega'$ with $\abs{\Omega'}=\omega_0$ cannot be translated into the reference domain $B$. 

In the first case, $\uei{1}$ is a minimizer of the functional $\Iei{0}$ and $\Omega(\uei{1})$ is an optimal domain for minimizing the fundamental tone in the sense of Theorem \ref{theo:main}. Thus, in this case, we can treat the functionals $\Iei{0}$ and $\Iei{1}$ as equivalent. 

In the second case, as a worst case scenario, we may think of the domain $\Omega(\uei{1})$ as of a domain with thin tentacles, which may all touch the boundary of the reference domain $B$. 
These tentacles eludes the scaling. Consequently, in this case a more local analysis of $\partial\Omega(\uei{1})$ is needed. 
Exemplary, we will see that assuming that $\partial\Omega(\uei{1})$ satisfies a doubling condition, the domain $\Omega(\uei{1})$ fulfills $\abs{\Omega(\uei{1})}=\omega_0$.

We begin this section with the analog result to Theorem \ref{theo:vol_0_2}.

\begin{theorem}\label{theo:vol_1_1}
  For $\eps\leq \eps_1$ every minimizer $\uei{1}$ of $\Iei{1}$ satisfies $\abs{\Omega(\uei{1})}\leq \omega_0$.  Thereby, $\eps_1$ is the number given in Theorem \ref{theo:vol_0_2}. 
\end{theorem}
\begin{proof}
  Let us assume that $\eps\leq \eps_1$ and that $\abs{\Omega(\uei{1})}=\alpha \omega_0$ for an $\alpha >1$. Then arguing in exactly the same way as in the proof of Theorem \ref{theo:vol_0_2} leads to a contradiction.  
\end{proof}

 We would like to repeat the above proof assuming that $\abs{\Omega(\uei{1})}=\alpha\omega_0$ for some $\alpha\in(0,1)$ and then define the scaled domain $\Omega' = \alpha^{-\frac{1}{n}}\Omega(\uei{1})$. However, it is not clear if the enlarged domain $\Omega'$ is still contained in $B$. A partial result can be obtained with the help of an argument of M. S. Ashbaugh and R. S. Laugesen. 

\begin{remark}\label{rem:AshLau}
  In \cite{AshLau1996}, M. S. Ashbaugh and R. S. Laugesen showed that there exists a constant $d_n \in (0,1)$, only depending on the dimension $n$, such that for every domain $\Omega \subset \R^n$ there holds 
\[
  \G(\Omega) > d_n \, \G(\Omega^\#),
\]
where $\Omega^\#$ denotes a ball in $\R^n$ with the same volume than $\Omega$. In addition, $d_n$ tends to $1$ as $n$ tends to infinity. 
\end{remark}

\begin{theorem}\label{theo:vol_1_2}
   Let $\uei{1}$ be a minimizer of $\Iei{1}$. Then there exists a number $\alpha_0 = \alpha_0(n,\eps_1,\eps)$ such that 
   \[
     \abs{\Omega(\uei{1})} \geq \alpha_0 \omega_0. 
   \] 
   Moreover, we have the explicit representation
   \[
      \alpha_0 = \frac{1+\eps\eps_1-\sqrt{1+2\eps\eps_1+(\eps\eps_1)^2-4d_n\eps\eps_1}}{2\eps\eps_1}, 
   \]
   where $\eps_1$ is given by Theorem \ref{theo:vol_0_2} and $d_n$ is given by Remark \ref{rem:AshLau}.
\end{theorem}

Note that 
\[
  \lim_{n\to\infty}\alpha_0 = 1 \quad \mbox{and } \lim_{\eps\to0}\alpha_0 = d_n.
\]
Together with Theorem \ref{theo:vol_1_1} this shows that for $\eps\leq\eps_1$ the domain $\Omega(\uei{1})$ satisfies the volume condition asymptotically as the dimension $n$ approaches infinity. 

\begin{proof}[Proof of Theorem \ref{theo:vol_1_2}]
 Let $\uei{1}$ be a minimizer of $\Iei{1}$ and let $\alpha \in (0,1)$ be such that 
 \[
    \abs{\Omega(\uei{1})} = \alpha \omega_0. 
 \]
 Choose $\phi$ as in Remark \ref{rem:EF_omega_0}. Then the minimality of $\uei{1}$ for $\Iei{1}$ implies
 \begin{equation}\label{eq:low_1}
   \G(\Omega(\uei{1})) - \eps\omega_0 (1-\alpha) = \Iei{1}(\uei{1}) \leq \Iei{1}(\phi) = \left(\frac{\omega_n}{\omega_0}\right)^\frac{4}{n}\G(B_1). 
 \end{equation}
 By $\Omega(\uei{1})^\#$ we denote the ball centered in the origin, having the same volume as $\Omega(\uei{1})$. Then 
 \[
   \G(\Omega(\uei{1})^\#) = \left( \frac{\omega_n}{\alpha\omega_0}\right)^\frac{4}{n}\G(B_1)
 \]
 and applying Remark \ref{rem:AshLau} we obtain from \eqref{eq:low_1}
 \begin{align*}
   &\left( \frac{\omega_n}{\omega_0}\right)^\frac{4}{n}\G(B_1) \left(d_n\alpha^{-\frac{4}{n}}-1\right) < \eps\,\omega_0\,(1-\alpha) \\
   \Rightarrow \quad &\frac{d_n\alpha^{-\frac{4}{n}}-1}{1-\alpha} < \eps\,\left(\frac{\omega_0}{\omega_n}\right)^\frac{4}{n}\frac{\omega_0}{\G(B_1)} = \eps\,\eps_1.
 \end{align*}
 Since $\alpha<1$, there holds $\alpha^{-\frac{4}{n}}\geq \alpha^{-1}$ and we have 
 \begin{equation}\label{eq:low_2}
   \frac{d_n\alpha^{-1}-1}{1-\alpha} \leq \eps\,\eps_1. 
 \end{equation}
 We set $f(\alpha) :=  \frac{d_n\alpha^{-1}-1}{1-\alpha}$. Then $f: (0,1)\to\R$ is smooth, strictly decreasing and 
 \[
   \lim_{\alpha \to 1}f(\alpha) = -\infty \; \mbox{and } \lim_{\alpha\to0}f(\alpha)=\infty. 
 \]
 By $\alpha_0$ we denote the (unique) solution in $(0,1)$ for equality in \eqref{eq:low_2}, i.e. 
 \[
      \alpha_0 = \frac{1+\eps\eps_1-\sqrt{1+2\eps\eps_1+(\eps\eps_1)^2-4d_n\eps\eps_1}}{2\eps\eps_1}.
   \]
Consequently, the strict monotonicity of $f$ implies that \eqref{eq:low_2} can only hold true for $\alpha \in [\alpha_0,1)$. This proves the theorem. 
\end{proof}

From now on, we always consider $0<\eps\leq\eps_1$. Consequently, there holds 
\[
   \alpha_0\omega_0 \leq \abs{\Omega(\uei{1})} \leq\omega_0.
\] 

\begin{theorem}\label{theo:vol_dicho}
  There exists a number $\eps_0 = \eps_0(n,\omega_0)$ such that for $\eps\leq \eps_0$ every minimizer $\uei{1}$ of $\Iei{1}$ satisfies either 
 \begin{enumerate}[label=\alph*)]
    \item $\abs{\Omega(\uei{1})} = \omega_0$  \qquad or 
    \item $\abs{\Omega(\uei{1})}< \omega_0$ and the rescaled domain $t\,\Omega(\uei{1})$ with $\abs{t\,\Omega(\uei{1})}=\omega_0$ is not a subset of $B$. In addition, there exists no translation $\Phi : \R^n\to\R^n$ such that $\Phi(t\,\Omega(\uei{1}))$ is contained in $B$.
   \end{enumerate}
\end{theorem}
\begin{proof}
  We claim that the statement of the theorem holds true for 
  \begin{equation}\label{eq:eps_0}
    \eps_0 := \min\left\{\eps_1, d_n\frac{\G(B_1)}{\omega_0}\frac{4}{n}\left(\frac{\omega_n}{\omega_0}\right)^\frac{4}{n}\right\} = \min\left\{ \eps_1, d_n\frac{4}{n}\eps_1^{-1}\right\},
  \end{equation} 
  where $\eps_1$ is given in Theorem \ref{theo:vol_0_2}  and $d_n$ is given in Remark \ref{rem:AshLau}. 
  We prove by contradiction.   Let $0<\eps\leq\eps_0$ and $\uei{1}$ be a minimizer of $\Iei{1}$.  Recall that Theorem \ref{theo:vol_1_2} implies $\abs{\Omega(\uei{1})} \in [\alpha_0\omega_0,\omega_0]$ since we assume $\eps\leq\eps_1$. Now let us assume that there exists an $\alpha \in [\alpha_0,1)$ such that 
   \begin{equation}\label{eq:dicho_0}
      \abs{\Omega(\uei{1})} = \alpha\omega_0 . 
   \end{equation}
  We split the proof in two steps.   
  
 Step 1. \quad Assume $\eps \leq \eps_0$ and let $\uei{0}$ be a minimizer of $\Iei{0}$. Recall that Theorem \ref{theo:vol_0_2} implies that $\abs{\Omega(\uei{0})}=\omega_0$. We now show that, for every $t \in (0,1)$, the clamped eigenfunction $u_{\eps}^t$  on the scaled domain $t\,\Omega(\uei{0})$ does not minimize the functional $\Iei{1}$, i.e. for every $t \in (0,1)$ there holds 
  \begin{equation}\label{eq:dicho_1}
     \Iei{1}(\uei{1}) < \Iei{1}(u_\eps^t).  
  \end{equation}
  For that purpose, let us assume that \eqref{eq:dicho_1} does not hold true. Then there exists a $t\in(0,1)$ such that 
  \[
    \Iei{1}(\uei{1}) = \Iei{1}(u_\eps^t)
  \] 
  and, since $\abs{t\,\Omega(\uei{0})}=t^n\omega_0 < \omega_0$, we obtain
  \begin{equation}\label{eq:dicho_2}
     \Iei{1}(u_\eps^t) = \G(t\,\Omega(\uei{0})) - \eps\omega_0(1-t^n) \leq \G(\Omega(\uei{0})) = \Iei{1}(\uei{0}).
  \end{equation}
  By scaling there holds 
  \[
    \G(t\,\Omega(\uei{0})) = t^{-4} \G(\Omega(\uei{0})).
  \]
  Thus, from \eqref{eq:dicho_2} we get 
  \begin{equation}\label{eq:dicho_2b}
    \G(\Omega(\uei{0})) \left(t^{-4} - 1\right) \leq \eps\omega_0 (1-t^n). 
  \end{equation}
  By $\Omega(\uei{0})^\#$ we denote the ball centered in the origin with the same volume as $\Omega(\uei{0})$. Then 
  \[
    \G(\Omega(\uei{0})^\#) = \left(\frac{\omega_n}{\omega_0}\right)^\frac{4}{n}\G(B_1)
  \]
  and applying Remark \ref{rem:AshLau} we obtain from \eqref{eq:dicho_2b}
  \begin{align*}
     &d_n \G(\Omega(\uei{0})^\#) \left(t^{-4}-1\right) < \eps \omega_0 \left(1-t^n\right)  \\
      \aq  \quad &d_n\,\left(\frac{\omega_n}{\omega_0}\right)^\frac{4}{n}\G(B_1)\left(t^{-4}-1\right) < \eps\,\omega_0\,(1-t^n) \\
   \aq \quad &\frac{t^{-4}-1}{1-t^n} < \eps \, \frac{\omega_0}{\G(B_1)}\,d_n^{-1}\,\left( \frac{\omega_0}{\omega_n}\right)^\frac{4}{n} = \eps\cdot\eps_1\cdot d_n^{-1}.
 \end{align*}
 Since $t <1$ and $\eps\leq\eps_0$ (see \eqref{eq:eps_0}),  we obtain
 \[
   \frac{4}{n} \leq \frac{t^{-4}-1}{1-t^n} < \eps\cdot\eps_1\cdot d_n^{-1} \leq \frac{4}{n}.
 \]
Obviously, this statement is false and we conclude that \eqref{eq:dicho_1} holds true for every $t \in (0,1)$.

Step 2. \quad Let us fix $t_\ast \in (0,1)$ such that $\abs{t_\ast\Omega(\uei{0})} = \abs{\Omega(\uei{1})}$. Hence, according to \eqref{eq:dicho_0}, there holds $t_\ast = \alpha^\frac{1}{n}$ and applying \eqref{eq:dicho_1} we obtain
\[
    \Iei{1}(\uei{1}) < \Iei{1}(u_\eps^{t_\ast}).
\]
By choice of $t_\ast$, this is equivalent to 
\[
  \G(\Omega(\uei{1})) < \G(t_\ast\Omega(\uei{0})) = t_\ast^{-4} \G(\Omega(\uei{0})),
\]
where we used the scaling property of the fundamental tone. 
Multiplying the above inequality with $t_\ast^4$ and again applying the scaling property of the fundamental tone, we deduce
\begin{equation}\label{eq:dicho_3}
   \G(t_\ast^{-1}\Omega(\uei{1})) = t_\ast^{-4}\G(\Omega(\uei{1})) < \G(\Omega(\uei{0})).
\end{equation}
Note that $\abs{t_\ast^{-1}\Omega(\uei{1})}=\omega_0$. Thus, if $t_\ast^{-1}\Omega(\uei{1}) \subset B$, \eqref{eq:dicho_3} is contradictory to the minimality of $\Omega(\uei{0})$ for the  fundamental tone among all open subsets of $B$ with volume less or equal than $\omega_0$ (see \eqref{eq:opt1}). Consequently, if $t_\ast^{-1}\Omega(\uei{1}) \subset B$, the assumption \eqref{eq:dicho_0} is false and there holds $\abs{\Omega(\uei{1})}=\omega_0$. This proves part a) of the claim. 

If $t_\ast^{-1}\Omega(\uei{1}) \not\subset B$, but there exists a translation $\Phi: \R^n\to\R^n$ such that $\Phi(t_\ast^{-1}\Omega(\uei{1})) \subset B$, we arrive at the same contradiction to \eqref{eq:opt1} as above because of the translational invariance of the fundamental tone (i.e. $\G(D) = \G(\Phi(D))$ for every  translation $\Phi$). 

Hence, if $\abs{\Omega(\uei{1})}<\omega_0$, the scaled domain $t_\ast^{-1}\Omega(\uei{1})$ cannot be translated into the ball $B$. This proves part b) of the claim.
\end{proof}

 From now on, we always choose $\eps\leq\eps_0$. Let $\uei{1}\in H^{2,2}_0(B)$ be a minimizer of $\Iei{1}$ and let us assume that case a) of Theorem \ref{theo:vol_dicho} holds true. Hence, $\abs{\Omega(\uei{1})}=\omega_0$ and repeating the proof of Theorem \ref{theo:main} we  gave in Section \ref{sec:nonrew} we find that $\Omega(\uei{1})$ minimizes the  fundamental tone among all open subsets of $B$ with volume smaller or equal than $\omega_0$. In addition, $\uei{1}$ minimizes the functional $\Iei{0}$.
 
 Consequently, if we could exclude that the case b) of Theorem \ref{theo:vol_dicho} may occur, we could treat the functionals $\Iei{0}$ and $\Iei{1}$ as equivalent.  We will discuss this issue in the following section.

\subsection{Discussion of case b) of Theorem \ref{theo:vol_dicho}}\label{sec:case_b}

From now on, we always assume that  $\partial\Omega(\uei{1})$ satisfies the following doubling property. Assume that there exists a constant $\sigma>0$ and a radius $0<R_0<1$ such that  for every $x_0 \in \partial\Omega(\uei{1})$ and every $0<R\leq R_0$ there holds 
\begin{equation}\label{eq:doubprop}
  \abs{B_{2R}(x_0)\cap\Omega(\uei{1})} \leq \sigma \, \abs{B_R(x_0) \cap \Omega(\uei{1})}.
\end{equation}
Note carefully that the condition \eqref{eq:doubprop} does not exclude that $\Omega(\uei{1})$ forms thin tentacles but it determines the minimal rate at which the volume of a possible tentacle may decrease.

\subsubsection{Nondegeneracy of $\uei{1}$}
Our aim is to establish a nondegeneracy result for $\uei{1}$.
For convenience, we cite a technical lemma which will be applied in the proof of Lemma~\ref{la:nondeg}. 

\begin{lemma}[{\cite[Lemma 6.1]{Giusti_direct_methods}}]\label{la:giusti}
Let $Z(t)$ be a bounded nonnegative function in the interval $[\rho,R]$. Assume that for $\rho\leq t < s \leq R$ we have
\[
 Z(t) \leq \left[A(s-t)^{-\alpha} + B(s-t)^{-\beta}+C\right] + \vartheta Z(s)
\]
with $A,B,C \geq 0$, $\alpha>\beta>0$ and $0\leq\vartheta<1$. Then,
\[
  Z(\rho) \leq c(\alpha,\theta)\left[ A(R-\rho)^{-\alpha} + B(R-\rho)^{-\beta}+C \right].
\]
\end{lemma}

\begin{lemma}\label{la:nondeg}
 Let $\uei{1}$ be a minimizer of $\Iei{1}$ and  let $\partial\Omega(\uei{1})$ satisfy \eqref{eq:doubprop}.
  There exists a $c_1 = c_1(\eps,n,\omega_0,\sigma) >0$  such that there holds 
  \begin{equation}\label{eq:nondeg_claim}
     c_1 \,R \leq \sup_{B_R(x_0)}\abs{\grad\uei{1}}, 
  \end{equation}
  where $x_0 \in \partial\Omega(\uei{1})$ and $0<R\leq R_0$.
  In particular, $c_1$ is independent of the choice of $x_0$ and $R$. 
\end{lemma} 
\begin{proof}
  Let $x_0 \in \partial\Omega(\uei{1})$ and $0<R\leq R_0$. 
  For $\frac{R}{2} \leq t < s \leq R$, let  $\eta \in C^\infty(\R^n)$ satisfy $0 \leq \eta \leq 1$, $\eta \equiv 0$ in $B_t(x_0)$ and $\eta \equiv 1$ in $R^n\setminus B_s(x_0) $.  Note that for every $x \in B_s(x_0)\setminus B_t(x_0)$ there holds 
  \begin{equation}\label{eq:eta}
     \abs{\eta(x)} \leq \frac{C(n)}{s-t} \; \mbox{and } \abs{\lap\eta(x)}\leq \frac{C(n)}{(s-t)^2}, 
  \end{equation}
  where $C(n)$ only depends on $n$.  We use $\uei{1}\eta$ as  a comparison function for $\Iei{1}$. By the minimality of $\uei{1}$ we obtain
  \[
    \Iei{1}(\uei{1}) \leq \Iei{1}(\uei{1}\eta). 
  \]
   Since, by construction, there holds  $\abs{\O(\eta\uei{1})}= \abs{\O(\uei{1})} - \abs{B_t(x_0)\cap \Omega(\uei{1})}$, straight forward computation yields
  \begin{equation}\label{eq:inner_cond_1}\begin{split}
    \int\limits_{B_s(x_0)}\abs{\lap \uei{1}}^2dx + \eps\,\abs{B_t&(x_0)\cap\Omega(\uei{1})} \leq \int\limits_{B_s\setminus B_t(x_0)}\abs{\lap(\eta\uei{1})}^2dx \\ &+ (\G(\Omega(\uei{1})+\eps\abs{B_t(x_0)\cap\Omega(\uei{1})})\int\limits_{B_s(x_0)} \uei{1}^2dx. 
  \end{split}\end{equation}
  A detailed analysis of this inequality gives the claim. This is done in three steps. 
  
  Step 1. \quad We estimate the last summand on the right hand side of \eqref{eq:inner_cond_1}. Recall that 
  \begin{align*}
     \G(\Omega(\uei{1})) &\leq \left(\frac{\omega_n}{\omega_0}\right)^\frac{4}{n}\G(B_1) + \eps\,(\omega_0 - \abs{\O(\uei{1})})\\
     &\leq \left(\frac{\omega_n}{\omega_0}\right)^\frac{4}{n}\G(B_1) + \eps_1\,\omega_0 = C(n,\omega_0). 
  \end{align*}
  Consequently, 
  \[
    \G(\Omega(\uei{1})) + \eps\abs{B_t(x_0)\cap\Omega(\uei{1})} \leq C(n,\omega_0) + \eps_1\abs{B_{R_0}(x_0)} \leq C(n,\omega_0)
  \]
  and \eqref{eq:inner_cond_1} becomes
  \begin{equation}\label{eq:inner_cond_2}
  \begin{split}
       \int\limits_{B_s(x_0)}\abs{\lap \uei{1}}^2dx + \eps\,\abs{B_t(x_0)&\cap\Omega(\uei{1})} \\ &\leq \int\limits_{B_s\setminus B_t(x_0)}\abs{\lap(\eta\uei{1})}^2dx + C(n,\omega_0)\int\limits_{B_s(x_0)} \uei{1}^2dx. 
  \end{split}
  \end{equation}
  
  Step 2. \quad Applying Young's inequality we estimate 
  \[
    \abs{\lap(\uei{1}\eta)}^2 \leq 4\left(\abs{\lap \uei{1}}^2\eta^2 + 2\abs{\grad\uei{1}.\grad\eta}^2 + \uei{1}^2\abs{\lap\eta}^2\right)
  \]
  and together with \eqref{eq:eta} we deduce in $B_s(x_0)\setminus B_t(x_0)$
  \[
     \abs{\lap(\uei{1}\eta)}^2 \leq 4\left(\abs{\lap\uei{1}}^2\eta^2 + 2\frac{C(n)}{(s-t)^2}\abs{\grad\uei{1}}^2 + \frac{C(n)}{(s-t)^4}\uei{1}^2\right). 
  \]
  With the splitting 
  \[
     \int\limits_{B_s(x_0)}\abs{\lap\uei{1}}^2dx = \int\limits_{B_s\setminus B_t(x_0) }\abs{\lap\uei{1}}^2  dx + \int\limits_{B_t(x_0)}\abs{\lap\uei{1}}^2 dx
  \]
  \eqref{eq:inner_cond_2} becomes
  \begin{align*}
     &\int\limits_{B_t(x_0)}\abs{\lap\uei{1}}^2dx + \eps\,\abs{B_t(x_0)\cap\Omega(\uei{1})} \leq \int\limits_{B_s\setminus B_t(x_0)}\abs{\lap\uei{1}}^2(4\eta^2-1)dx \\
     &+ \int\limits_{B_s\setminus B_t(x_0)}\frac{C(n)}{(s-t)^2}\abs{\grad\uei{1}}^2 + \frac{C(n)}{(s-t)^4}\uei{1}^2dx + C(n,\omega_0)\int\limits_{B_s(x_0)}\uei{1}^2dx,
  \end{align*}
  where $C(n)$ collects all constants only depending on $n$. Since $0\leq\eta\leq 1$, we obtain
\begin{align*}
     &\int\limits_{B_t(x_0)}\abs{\lap\uei{1}}^2dx + \eps\,\abs{B_t(x_0)\cap\Omega(\uei{1})} \leq 3 \int\limits_{B_s\setminus B_t(x_0)}\abs{\lap\uei{1}}^2dx \\
     &\quad+C(n,\omega_0)\, \int\limits_{B_R(x_0)}\frac{1}{(s-t)^2}\abs{\grad\uei{1}}^2 +\left(1+ \frac{1}{(s-t)^4}\right)\uei{1}^2dx .
\end{align*}
Now we add $3\int_{B_t(x_0)}\abs{\lap\uei{1}}^2dx$ to both sides of the above inequality and divide the resulting inequality by $4$. Subsequently, we add $\frac{3\,\eps}{16}\abs{B_s(x_0)\cap\Omega(\uei{1})}$ to the right hand side. This leads to
\begin{equation}\begin{split}\label{eq:inner_cond_3}
    &\int\limits_{B_t(x_0)}\abs{\lap\uei{1}}^2dx + \frac{\eps}{4}\,\abs{B_t(x_0)\cap\Omega(\uei{1})} \leq \frac{3}{4}\left( \int\limits_{B_s(x_0)}\abs{\lap\uei{1}}^2dx + \frac{\eps}{4}\abs{B_s(x_0)\cap\Omega(\uei{1})}\right) 
    \\& \quad+C(n,\omega_0)\, \int\limits_{B_R(x_0)}\frac{1}{(s-t)^2}\abs{\grad\uei{1}}^2 + \left(1+ \frac{1}{(s-t)^4}\right)\uei{1}^2dx .
\end{split}\end{equation}
Setting 
\[
  Z(t) := \int\limits_{B_t(x_0)}\abs{\lap\uei{1}}^2dx + \frac{\eps}{4}\,\abs{B_t(x_0)\cap\Omega(\uei{1})},
\]
estimate \eqref{eq:inner_cond_3} enables us to apply Lemma \ref{la:giusti} and  we obtain
\[
     \int\limits_{B_\frac{R}{2}(x_0)}\abs{\lap\uei{1}}^2dx + \frac{\eps}{4}\,\abs{B_\frac{R}{2}(x_0)\cap\Omega(\uei{1})} \leq
   C(n,\omega_0)\, \int\limits_{B_R(x_0)}\frac{1}{R^2}\abs{\grad\uei{1}}^2 + \left(1+\frac{1}{R^4}\right)\uei{1}^2dx.
\]
Step 3.\quad The $C^{1,\alpha}$ regularity of $\uei{1}$ allows us to estimate 
\[
   \abs{\uei{1}(x)} \leq 2\,R\,\sup_{B_R(x_0)}\abs{\grad \uei{1}} 
\]
for every $x\in B_R(x_0)$. Moreover, we assume $R\leq R_0<1$. Hence, we find
\begin{equation}\label{eq:why_dp}
\begin{split}
      \int\limits_{B_\frac{R}{2}(x_0)}\abs{\lap\uei{1}}^2dx + \frac{\eps}{4}\,\abs{B_\frac{R}{2}(x_0)&\cap\Omega(\uei{1})} \\ & \leq
   C(n,\omega_0)\, R^{-2}\,\sup_{B_R(x_0)}\abs{\grad \uei{1}}^2\,\abs{B_R(x_0)\cap\Omega(\uei{1})}.
\end{split}
\end{equation}
Since we assume the doubling property \eqref{eq:doubprop} to hold true, omitting the nonnegative integral on the left hand side we obtain
\[
 \frac{\eps}{4} \leq C(n,\omega_0)\,\sigma\, R^{-2}\,\sup_{B_R(x_0)}\abs{\grad \uei{1}}^2.
\]
This proves the claim.
\end{proof}

Note that the assumption \eqref{eq:doubprop} enables us to compare $\abs{B_\frac{R}{2}(x_0)\cap\Omega(\uei{1})}$ with $\abs{B_R(x_0)\cap\Omega(\uei{1})}$ in estimate \eqref{eq:why_dp}. This assumption is only needed because we currently do not have any further information about the free boundary $\partial \Omega(\uei{1})$ and could be replaced by regularity properties of the free boundary. However, let us emphasize that the rewarding property of the penalization term $\pen{1}$ is crucial for proving Lemma \ref{la:nondeg} and cannot be replaced since the rewarding term yields the strictly positive lower bound in \eqref{eq:nondeg_claim}. 

The nondegeneracy of $\uei{1}$ along the free boundary according to Lemma \ref{la:nondeg} allows us to establish a lower bound on the density quotient of $\Omega(\uei{1})$. 

\begin{lemma}\label{la:density}
  Let $\uei{1}$ be a minimizer of $\Iei{1}$, $\alpha \in (0,1)$ and let $\partial\Omega(\uei{1})$ satisfy \eqref{eq:doubprop}. 
  There exists a constant $c_2 = c_2(\eps,n,\omega_0,\sigma,\alpha)>0$ such that for every $x_0 \in \partial\Omega(\uei{1})$ and every $0<R\leq R_0$ there holds 
  \[
     c_2 \,\abs{B_R}^\frac{1-\alpha}{\alpha} \leq \frac{\abs{\Omega(\uei{1})\cap B_R(x_0)}}{\abs{B_R}}.
  \]
\end{lemma} 
Although this lower bound on the density quotient is admittedly weak, it suffices to prove that $\Omega(\uei{1})$ can be rescaled to the volume $\omega_0$ without leaving the reference domain $B$ provided the radius of $B$ is chosen sufficiently large (see Theorem \ref{theo:radius}). 

\begin{proof}[Proof of Lemma \ref{la:density}.] 
  Let $x_0 \in \partial\Omega(\uei{1})$ and $0<R\leq R_0$. According to Lemma \ref{la:nondeg} and  the $C^{1,\alpha}$ regularity of $\uei{1}$ there exists an $x_1 \in \Omega(\uei{1})\cap\overline{B_\frac{R}{2}(x_0)}$ such that 
 \[
    c_1 \frac{R}{2} \leq \sup_{B_\frac{R}{2}(x_0)}\abs{\grad\uei{1}} = \abs{\grad \uei{1}(x_1)}.
 \]
 Now choose $x_2 \in \partial\Omega(\uei{1})$ such that 
 \[
   d:=  \dist(x_1,\partial\Omega(\uei{1}) =  \abs{x_1-x_2}.
 \]
 Since $\abs{\grad \uei{1}(x_2)}=0$, we obtain 
 \[
    c_1 \frac{R}{2} \leq \abs{\grad\uei{1}(x_1)-\grad\uei{1}(x_2)} \leq L_\alpha\,d^\alpha,
 \]
 where $L_\alpha=L(n,\omega_0,\alpha)$ denotes the $\alpha$-Hölder coefficient of $\grad\uei{1}$. By construction, there holds  $B_d(x_1)\subset \Omega(\uei{1})\cap B_R(x_0)$. Consequently, we may proceed to
 \begin{align*}
     \left(\frac{c_1}{2L_\alpha}\right)^n\,R^n \leq d^{\alpha\,n} &\aq \left(\frac{c_1}{2L_\alpha \omega_n^{1-\alpha}}\right)^n\,\abs{B_R} \leq \abs{B_d(x_1)}^\alpha \\
     &\Rightarrow  \left(\frac{c_1}{2L_\alpha \omega_n^{1-\alpha}}\right)^n\,\abs{B_R} \leq \abs{\Omega(\uei{1})\cap B_R(x_0)}^\alpha.
 \end{align*}
 This proves the claim.
\end{proof}

\subsubsection{The volume condition for $\Omega(\uei{1})$}

The next theorem is the key observation to show that $\Omega(\uei{1})$ has the volume $\omega_0$ provided that $\partial\Omega(\uei{1})$ satisfies \eqref{eq:doubprop}.  
It is a consequence of the lower bound on the density quotient according to Lemma \ref{la:density}. 

\begin{theorem}\label{theo:radius}
 Let $\eps\leq\eps_1$ and let $B = B_{R_B}(0)$. Provided that $R_B$ is chosen sufficiently large, for every minimizer $\uei{1}$ of $\Iei{1}$ such that $\partial\Omega(\uei{1})$ satisfies \eqref{eq:doubprop} there  holds 
\begin{enumerate}[label=\alph*)]
  \item $\Omega(\uei{1})$ is compactly contained in $B_{2^{-\frac{1}{n}}R_B}(0)$ \; or
  \item there exists a translation $\Phi: \R^n\to\R^n$ such that $\Phi(\Omega(\uei{1}))$ is compactly contained in $B_{2^{-\frac{1}{n}}R_B}(0)$. 
\end{enumerate}
\end{theorem}
\begin{proof}
Let us think of the reference domain $B$ as of a ball centered at the origin and with radius $R_B$.
In addition, let $\uei{1}$ be a minimizer of $\Iei{1}$ for $\eps\leq\eps_0$. 
In order to prove the claim, let us assume that $\Omega(\uei{1})$ is not compactly contained in $B_{2^{-\frac{1}{n}}R_B}(0)$. For the sake of convenience, we abbreviate $S := 2^{-\frac{1}{n}}R_B$. Hence, we assume $\partial\Omega(\uei{1})\cap\partial B_{S}(0) \neq \emptyset$. 

Of course, there either holds $0 \in \Omega(\uei{1})$ or $0\not\in\Omega(\uei{1})$. At first, we will handle the case where the origin is already contained in $\Omega(\uei{1})$. Secondly, we will show that we may translate $\Omega(\uei{1})$ such that the origin becomes an inner point of $\Omega(\uei{1})$. 

Step 1. \quad We consider that $0\in\Omega(\uei{1})$. 
  
Note that for every $m\in\N$ with $m\geq 3$ there holds 
  \begin{equation}\label{eq:B_covering}
      B_S(0) = \bigcup_{i=0}^{m-2} B_{\frac{i+2}{m}S}(0) \setminus B_{\frac{i}{m}S}(0). 
  \end{equation}
  Since we assume that $\partial\Omega(\uei{1})\cap\partial B_S(0) \neq \emptyset$ and $0\in\Omega(\uei{1})$, there exists a smallest index $i_0=i_0(m)$ such that for each $i\geq i_0$ there exists an $x_i \in \partial\Omega(\uei{1})\cap\partial B_{\frac{i+1}{m}S}(0)$.  
  
 We fix $m\in\N$ such that $\frac{S}{m} \leq R_0.$ In addition, we fix an $\alpha \in (0,1)$.
 Then applying Lemma \ref{la:density}  for  $i_0\leq i\leq m-2$ we obtain
\begin{equation}\label{eq:radius_1}
c_2 \abs{B_\frac{S}{m}}^\frac{1}{\alpha} \leq \abs{\Omega(\uei{1})\cap B_\frac{S}{m}(x_i)}. 
\end{equation}
We now sum \eqref{eq:radius_1} from $i=i_0(m)$ to $i=m-2$. Since $B_\frac{S}{m}(x_i)\cap B_\frac{S}{m}(x_k)=\emptyset$ for $i\neq k$ and $\abs{\Omega(\uei{1})}\leq\omega_0$, this implies
\[
 c_2\,(m-1-i_0(m)) \, \abs{B_{\frac{S}{m}}}^\frac{1}{\alpha} \leq \sum_{i=i_0}^{m-2}\abs{\Omega(\uei{1})\cap B_\frac{S}{m}(x_i)} \leq \abs{\Omega(\uei{1})} \leq \omega_0.
\]
Note that since $B_{\frac{i_0}{m}S}(0)\subset\Omega(\uei{1})$ and $\abs{\Omega(\uei{1})}\leq\omega_0$,  $i_0(m)$ is bounded. 
Indeed, 
\[
\abs{B_{\frac{i_0}{m}S}} \leq \abs{\Omega(\uei{1})} \leq \omega_0 
\]
implies
\[
  i_0(m) \leq \left(\frac{\omega_0}{\omega_n}\right)^\frac{1}{n}\frac{m}{S}.
\]
Specifying the choice of $m\in \N$ such that 
\[
   \frac{R_0}{4} \leq \frac{S}{m}\leq \frac{R_0}{2}
\]
and recalling that $S:=2^{-\frac{1}{n}}R_B $, we obtain
\[
   c_2\,\left(\frac{2\,R_B}{2^\frac{1}{n}R_0}-1-\frac{4}{R_0}\left(\frac{\omega_0}{\omega_n}\right)^\frac{1}{n}\right) \abs{B_\frac{R_0}{4}}^\frac{1}{\alpha} \leq c_2 \,(m-1-i_0(m))\abs{B_\frac{S}{m}}^\frac{1}{\alpha} \leq \omega_0.
\]
Since this estimate is false if $R_B$ is chosen sufficiently large, the proof is finished provided that $0 \in \Omega(\uei{1})$. 
 
 Step 2. \quad Let us now assume that the origin is not contained in $\Omega(\uei{1})$. However, there exists an $x_0 \in \Omega(\uei{1})\cap B$. We now translate $\Omega(\uei{1})$ such that $x_0$ is translated to the origin, i.e. we consider 
 \[
  \Phi : \R^n \to \R^n,\, x \mapsto x-x_0. 
 \]
 We call $\Omega' :=\Phi(\Omega(\uei{1}))$ and $v_\eps(x) := \uei{1}(\Phi^{-1}(x))$. Thus, 
 \[
   \Omega' = \{x\in\R^n : v_\eps(x)\neq 0 \mbox{ or } (v_\eps(x)=0 \wedge \abs{\grad v_\eps(x)}>0)\} 
 \]
 and $v_\eps\in H^{2,2}_0(\Omega')$. Moreover, $\partial\Omega' = \Phi(\partial\Omega(\uei{1}))$. 
 
 Let us emphasize that, in general, $\Omega'$ may not be contained in $B$ and, thus, $v_\eps \not\in H^{2,2}_0(B)$.
 Note carefully that in Step 1 the minimality of $\uei{1}$ for $\Iei{1}$ is not used explicitly.  However, the minimality is necessary to establish Lemma \ref{la:nondeg} and, subsequently, Lemma \ref{la:density} and the application of Lemma \ref{la:density} leads to estimate \eqref{eq:radius_1}, which is the crucial observation in Step 1. 
 
 If $y_0 \in \partial\Omega'$ and $0<R\leq R_0$, then there exists a $z_0 \in \partial\Omega(\uei{1})$ such that $y_0 = \Phi(z_0)$. Lemma \ref{la:density} together with the translational invariance of the Lebesgue measure then imply
\begin{equation}\label{eq:radius_5}
    \abs{\Omega' \cap B_R(y_0)} = \abs{\Omega(\uei{1})\cap B_R(z_0)} \geq c_2 \abs{B_R}^\frac{1}{\alpha}.
\end{equation}
 Estimate \eqref{eq:radius_5} enables us to repeat the approach presented in Step 1.
Again we consider the segmentation \eqref{eq:B_covering} and assume that $\partial\Omega'\cap\partial B_S(0)$ is not empty. Then there exists a smallest index $i_0(m)$ such that for every $i_0\leq i \leq m-2$ there exists an $x_i \in \partial\Omega'\cap\partial B_{\frac{i+1}{m}S}(0)$. Now applying \eqref{eq:radius_5}  we obtain for $i_0\leq i \leq m-2$
\[
c_2\,\abs{B_\frac{S}{m}}^\frac{1}{\alpha} \leq \abs{\Omega' \cap B_\frac{S}{m}(x_i)}
\]
Since $\abs{\Omega'} = \abs{\Omega(\uei{1})} \leq \omega_0$ we may repeat the argumentation from Step 1 and obtain that $\Omega'=\Phi(\Omega(\uei{1}))$ is compactly contained in $B_{2^{-\frac{1}{n}}R_B}$ and, by construction, contains the origin. 
 \end{proof}
 
 As a direct consequence of Theorem \ref{theo:radius} we deduce that, if  $\partial\Omega(\uei{1})$ satisfies the doubling condition \eqref{eq:doubprop}, the domain $\Omega(\uei{1})$ satisfies $\abs{\Omega(\uei{1})}=\omega_0$.
 
 \begin{corollary}
   Let $\eps\leq\eps_0$ and let $\uei{1}\in H^{2,2}_0(B)$ minimize $\Iei{1}$. Provided that the radius $R_B$ of $B$ is chosen sufficiently large and that $\partial\Omega(\uei{1})$ satisfies the doubling property \eqref{eq:doubprop}, there holds 
   $\abs{\Omega(\uei{1})}=\omega_0$.
 \end{corollary}
 \begin{proof}
     Let $\eps\leq\eps_0$ and let  $\uei{1}\in H^{2,2}_0(B)$ be a minimizer of $\Iei{1}$ such that $\partial\Omega(\uei{1})$ satisfies \eqref{eq:doubprop}. In addition, we assume $\abs{\Omega(\uei{1})}=\alpha\omega_0$ for an $\alpha \in [\alpha_0,1)$. Note that since $d_n\geq \frac{1}{2}$ for every $n \in \N$, the quantity $\alpha_0$ given in Theorem \ref{theo:vol_1_2} satisfies
     \begin{equation}\label{eq:bound_alpha_0}
        \alpha_0\geq \frac{1+\eps_1^2-\sqrt{(1+\eps_1^2)^2-4d_n\eps_1^2}}{2\eps_1^2} \geq \frac{1+\eps_1^2-\sqrt{1+\eps_1^4}}{2\eps_1^2} \geq \frac{1}{2} 
     \end{equation}
     for every choice of $\omega_0$, which determines $\eps_1$. Applying Theorem \ref{theo:radius} the domain $\Omega(\uei{1})$ is compactly contained in $B_{2^{-\frac{1}{n}}R_B}$ (possibly after translation) and estimate \eqref{eq:bound_alpha_0} implies that the scaled domain $\alpha^{-\frac{1}{n}}\Omega(\uei{1})$ is a subset of $B_{R_B}$. This is a contradiction to Theorem \ref{theo:vol_dicho} since we assume that $\abs{\Omega(\uei{1})}<\omega_0$.
 \end{proof}

\noindent\textbf{Funding.} The author is funded by the Deutsche Forschungsgemeinschaft (DFG, German Research Foundation) - project number 396521072.


\end{document}